\documentclass[11pt]{article}

\usepackage{amsfonts, amssymb}
\usepackage{amsmath}
\usepackage{verbatim}
\usepackage{indentfirst}
\usepackage{color}

\begin{document}
\title{Metric $n$-Lie Algebras}
\date{}
\maketitle
\vspace{ -1.7cm}

\begin{center}
    Ruipu Bai\footnote{$^{,~2}$ Project partially supported by NSF(10871192) of
China.  Email address:
bairp1@yahoo.com.cn; wuwanqing8888@126.com; zhenhengl@usca.edu.}  ~~Wanqing Wu$^2$ \\
\vspace{1mm}
College of Mathematics and Computer Science\\
Hebei University, Baoding (071002), China\\
\vspace{3mm}
Zhenheng Li\\
\vspace{1mm}
Department of Mathematical Sciences\\
University of South Carolina Aiken\\
Aiken, SC 29801, USA\\
\end{center}

\begin{abstract}
We study the structure of a metric $n$-Lie algebra $\mathcal {G}$
over the complex field $\mathbb C$. Let $\mathcal {G}= \mathcal
S\oplus {\mathcal R}$ be the Levi decomposition, where $\mathcal R$
is the radical of $\mathcal {G}$ and $\mathcal S$ is a strong
semisimple subalgebra of $\mathcal {G}$. Denote by $m(\mathcal {G})$
the number of all minimal ideals of an indecomposable metric $n$-Lie
algebra and $\mathcal R^\bot$ the orthogonal complement of $R$. We
obtain the following results. As $\mathcal S$-modules, $\mathcal
R^{\bot}$ is isomorphic to the dual module of $\mathcal {G} /
\mathcal R.$ The dimension of the vector space spanned by all
nondegenerate invariant symmetric bilinear forms on $\mathcal {G}$
equals that of the vector space of certain linear transformations on
$\mathcal {G}$; this dimension is greater than or equal to
$m(\mathcal {G}) + 1$. The centralizer of $\mathcal R$ in $\mathcal
G$ equals the sum of all minimal ideals; it is the direct sum of
$\mathcal R^\bot$ and the center of $\mathcal {G}$. The sufficient
and necessary condition for $\mathcal {G}$ having no strong
semisimple ideals is that $\mathcal R^\bot \subseteq \mathcal R$.

\vspace{ 0.3cm}
\noindent {\bf 2010 Mathematics Subject Classification:} 17B05 17D99

\vspace{2mm}\noindent{\bf Keywords:} Metric $n$-Lie algebra,
Minimal ideal, Metric dimension, Levi decomposition.
\end{abstract}

\baselineskip 18pt

\vspace{2mm}\noindent{\bf 1. Introduction}

\vspace{2mm} In mathematical community, $n$-Lie algebras were introduced in 1985 (see \cite{F}). These algebras
appear in many fields in mathematics and mathematical physics. Motivated by problems of quark dynamics,
 Nambu \cite{N} introduced in 1973 an $n$-ary generalization of Hamiltonian dynamics by means of the $n$-ary Poisson bracket
$$
[f_1, \cdots, f_n]= \det\Big(\frac{\partial f_i}{\partial x_j}\Big).
\eqno(1.1)
$$
Based on Nambu's work, Takhtajin \cite{T} built systematically the foundation of the theory of $n$-Poisson or
Nambu-Poisson manifolds. Physicists noticed that the identity (1.1) satisfies the generalized Jacobi identity given by  (2.2) below.

A metric $n$-Lie algebra is an $n$-Lie algebra that admits a nondegenerate symmetric bilinear form which is invariant
under inner derivations. Figueroa-O'Farrill and Papadopoulos introduced $n$-Lie algebras in their study of the classification
of maximally supersymmetric type IIB supergravity backgrounds (see \cite{FP,FP2,FO1}). Bagger and Lambert \cite{BL1,BL2,BL3}
and Gustavsson \cite{G} on superconformal fields for multiple M2-branes show that the generalized Jacobi identity for $3$-Lie
 algebras is essential in defining the action with N = 8 supersymmetry, and that this identity may also be considered a
 generalized Plucker relation in physics. Their work stimulates the interest of researchers in mathematics and mathematical physics on metric $n$-Lie algebras.

It is useful to study metric $n$-Lie algebras in both physical and mathematical observations. The ordinary gauge theory
 requires a positive-definite metric to guarantee that the theory possess positive-definite kinetic terms and to prevent
 violations of unitarity due to propagating ghost-like degrees of freedom. But very few metric $n$-Lie algebras admit
 positive-definite metrics (see \cite{P, JLZ}); Ho-Hou-Matsuo in \cite{HHM} confirmed that there are no $n$-Lie algebras
 with positive-definite metrics for n = 5, 6, 7, 8. They also gave examples of 3-Lie algebras whose metrics are not
 positive-definite and observed that generators of zero norm are common in 3-Lie algebras. Furthermore, they conjectured
 that there are no other 3-Lie algebras with a positive definite metric except for $\mathcal A_4$ and its direct sum, where
 $\mathcal A_4 $ is the $SO$(4)-invariant algebra with 4 generators. In physics, some dynamical systems involve zero-norm generators
 corresponding to gauge symmetries and negative-norm generators corresponding to ghosts (see \cite{MFO1,MFO2, MFO3}). These messages
  motivate us to study $n$-Lie algebras with any metrics or a nondegenerate invariant symmetric bilinear form. Note that Figueroa-O'Farrill \cite{FO2} classifies $n$-Lie
  algebras possessing an invariant lorentzian inner product.

We are interested in structures of metric $n$-Lie algebras $\mathcal {G}$ over the complex field. Section $2$  introduces  necessary notation and basic facts. Section
3 describes relations between metric dimensions and minimal ideals of $\mathcal {G}$. Section $4$ is devoted to study the structure of minimal ideals of metric $n$-Lie algebras.
 Throughout this paper, all $n$-Lie algebras are of finite dimension
and over the complex field $\mathbb C$.

\vspace{2mm}\noindent{\bf 2. Fundamental notions}

An {\it $n$-Lie algebra} is a vector space $\mathcal {G}$ over $\mathbb C$ equipped
 with an $n$-multilinear operation $[x_1, \cdots, x_n]$ satisfying
$$ [x_1, \cdots, x_n] = sign(\sigma)[x_{\sigma (1)}, \cdots, x_{\sigma(n)}], \eqno(2.1) $$  and
 $$
  [[x_1, \cdots, x_n], y_2, \cdots, y_n]=\sum_{i=1}^n[x_1, \cdots, [ x_i, y_2, \cdots, y_n], \cdots, x_n] \eqno(2.2)
$$
for $x_1, \cdots, x_n, y_2, \cdots, y_n\in \mathcal {G}$ and $\sigma\in S_n$, the permutation group on $n$ letters.  Identity (2.2) is usually called the
generalized Jacobi identity, or simply the Jacobi identity. The
subalgebra generated by all vectors $[x_1, \cdots, x_n]$ for $x_1,
\cdots, x_n\in \mathcal {G}$ is called the {\it derived algebra} of
$\mathcal {G}$, denoted by $\mathcal {G}^1$. If $\mathcal {G}^1 = \mathcal {G}$ then $\mathcal {G}$ is referred to as a perfect $n$-Lie algebra.

 An {\it ideal} of an $n$-Lie algebra $\mathcal {G}$ is a subspace $I$ such that $ [I,
\mathcal {G}, \cdots, \mathcal {G}]\subseteq I.$  If $I$ satisfies
$[I, I, \mathcal {G}, \cdots, \mathcal {G}]=0,$ then $I$ is called an {\it abelian ideal.} If $\mathcal {G}^1\neq 0$ and $\mathcal {G}$ has  no ideals except for $0$ and itself, then $\mathcal {G}$ is referred to as a {\it  simple $n$-Lie algebra}. By Ling \cite{L}, up to isomorphism there exists only one finite dimensional simple
 $n$-Lie algebra over an algebraically closed field of characteristic $0$; this algebra  is the $(n+1)$ dimensional $n$-Lie algebra whose derived algebra has dimension $n+1.$
 If an $n$-Lie algebra is the direct sum of its simple ideals,
 then it is called a {\it  strong semisimple $n$-Lie algebra} \cite{BM}. We call $\mathcal {G}$ indecomposable if $\mathcal {G}=J\oplus H$ is a direct sum of ideals then $J=\mathcal {G}$ or $H=\mathcal {G}$.

For a given subspace $V$ of an  $n$-Lie algebra $\mathcal {G}$, the subalgebra
$$
\mathcal C_{\mathcal {G}}(V)=\{x\in \mathcal {G} \mid [x, V,
\mathcal {G}, \cdots, \mathcal {G}]=0\}
$$
 is called the {\it  centralizer} of $V$ in $\mathcal {G}$. The centralizer of $\mathcal G$ in $\mathcal {G}$ is referred to
 as the {\it  center} of $\mathcal {G}$, and is denoted by $\mathcal C(\mathcal {G})$. It is clear that if $I$ is an ideal of $\mathcal {G}$, so is $\mathcal C_{\mathcal {G}}(I)$.

Let $I$ be an ideal of $\mathcal G$. Then $I$ is called a
characteristic ideal if $I$ is invariant under every derivation on $\mathcal G$. We call $I$ {\it solvable}, if $I^{(s)}=0$ for some
$s\geq 0$, where $I^{(0)}=I$ and $I^{(s)}$ is defined by induction,
$$
I^{(s+1)}=[I^{(s)},  \cdots, I^{(s)}]\eqno(2.3)
$$
for $s\geq 0$.  The maximal solvable ideal of $\mathcal {G}$ is
called the {\it  radical of $\mathcal {G}$}, denoted by $\mathcal R$.

A {\it minimal} ({\it maximal}) {\it ideal} is a non-trivial ideal
$I$ satisfying that if $J$ is an ideal of $\mathcal {G}$, $J\subset
I$ $(J\supset I)$ and  $J\ne I$ then $J=0$ ($J= \mathcal {G}$). The
$socle$ of $\mathcal {G}$ is the sum of all minimal ideals, denoted
by ${\rm Soc}(\mathcal {G})$. Then ${\rm Soc}(\mathcal {G})$=0 if
and only if $\mathcal {G}$ is simple or one-dimensional. An $n$-Lie
algebra $\mathcal {G}$ is local if it has only one minimal ideal.

  {\it A metric $n$-Lie
algebra } is an $n$-Lie algebra  $\mathcal G$ that possess a
nondegenerate symmetric bilinear form $B$ on $\mathcal G$, which is
invariant,
$$B([x_{1}, \cdots, x_{n-1}, y_{1}], y_{2})=-B([x_{1}, \cdots, x_{n-1}, y_{2}], y_{1}),
~\mbox{ for all} ~x_{i}, y_{j}\in \mathcal G.$$ Such a bilinear form
$B$ is called an {\it invariant scalar product} on $\mathcal {G}$ or
a {\it metric} on $ \mathcal {G}$. Note the $B$ is not necessarily
positive definite. Note pairs $(\mathcal {G}, B)$ be a {\it metric
$n$-Lie algebra}.

Let $W$ be a subspace of a metric $n$-Lie algebra $\mathcal G$. The orthogonal complement of $W$ is defined by
$$
    W^{\bot}=\{x \in \mathcal G \mid B(w, x)=0~ \mbox{for all} ~ w\in W\}.
$$
If $W$ is an ideal, then $W^{\bot}$ is also an ideal and $(
W^{\bot})^{\bot}=W$. Notice that $W$ is a minimal ideal if and only
if $W^{\bot}$ is maximal. We say that $W$ is {\it isotropic}
{\rm(}{\it coisotropic}{\rm)} if $W\subseteq W^{\bot}$  ($W^{\bot}
\subseteq W$). The subspace $W$ is called {\it nondegenerate} if
$B|_{W\times W}$ is nondegenerate; this is equivalent to  $W\cap
W^{\bot}=0$ or $ \mathcal G =W\oplus W^{\bot}$ as a direct sum of
subspaces. If $\mathcal {G}$ contains no nontrivial nondegenerate
ideals, then $\mathcal G$ is called {\it $B$-irreducible.} For a
metric $n$-Lie algebra $(\mathcal {G}, B)$, it is not difficult to
see $[\mathcal G, \cdots, \mathcal G]=\mathcal {C}(\mathcal
{G})^{\bot}$.

\vspace{2mm}\noindent{\bf Lemma 2.1}$^{\cite{L}}$  Let $\mathcal {G}$
be an $n$-Lie algebra. Then $\mathcal {G}$ has the Levi
decomposition
$$
    \mathcal {G}={\mathcal S}\oplus {\mathcal R} ~\mbox{ (as a direct sum of subspaces)}, \eqno(2.3)
$$
where $\mathcal R$ is the radical of $\mathcal {G}$ and $\mathcal S$ is a strong
semisimple subalgebra of $\mathcal {G}$.

\vspace{2mm}
Using induction on $\dim \mathcal G$, we have the following result without showing the details.

\vspace{1mm}\noindent{\bf Lemma 2.2 }  Every metric $n$-Lie algebra
$(\mathcal {G}, B)$ has a decomposition
$$\mathcal {G}=\oplus_{i=1}^{r}\mathcal G_{i} \quad\mbox{with}\quad [\mathcal G_{i}, \mathcal G_{j}, \mathcal G, \cdots, \mathcal G]=0 \quad\mbox{and}\quad B(\mathcal {G}_i, \mathcal {G}_j)=0 \quad\mbox{if} ~i\neq j,
$$
where $\mathcal G_{i}$ for $i=1, \cdots, r$ are $B$-irreducible nondegenerate ideals. Moreover, if $\mathcal {G}=\oplus_{k=1}^{l}\mathcal T_{k}$ is another decomposition of $B$-irreducible nondegenerate ideals, then $l=r$ and there exists a permutation $\sigma$ of $\{1, \cdots, r\}$ such that $\mathcal G_{i}$ is isomorphic to $\mathcal T_{\sigma(i)}$.

\vspace{2mm}\noindent{\bf Lemma 2.3 } Let $(\mathcal {G}, B)$ be a
metric $n$-Lie algebra and $I$ an ideal of $\mathcal {G}$. Then we
have following properties.

\vspace{2mm} (1) $[I, \mathcal {G}, \cdots, \mathcal
{G}]^{\bot}=\mathcal {C}_{\mathcal G}(I)$.

\vspace{2mm} (2) If $\mathcal G = I\oplus Y$ is a direct sum of ideals with $I=[I, \cdots, I]\neq 0 $ and $Y\neq 0$, then $I$ is nondegenerate. In particular, if $I$ is strong semisimple, then $I$ is nondegenerate.

 \vspace{2mm} (3)  If $I$ is nondegenerate, then  $\mathcal G=I\oplus I^{\bot}$.

\vspace{2mm} (4)  $\mathcal {G}$ is perfect if and only if $\mathcal
{C}(\mathcal {G})=0.$

\vspace{2mm}\noindent{\bf Proof } For every $y \in [I, \mathcal
{G}, \cdots, \mathcal {G}]^{\bot}$, we have $B([I, \mathcal {G},
\cdots, \mathcal {G}], y)= B(\mathcal {G}, [y, I, \mathcal {G},
\cdots, \mathcal {G}])=0$. Since $B$ is nondegenerate, $y\in \mathcal
{C}_{\mathcal G}(I)$. Conversely, for all $z\in \mathcal
{C}_{\mathcal G}(I)$, $[I, z, \mathcal {G}, \cdots, \mathcal
{G}])=0,$ and then $B(z, [I, \mathcal {G}, \cdots, \mathcal {G}])=
B(\mathcal {G},[I, z, \mathcal {G}, \cdots, \mathcal {G}])=0$. In other words, $z\in[I, \mathcal {G}, \cdots, \mathcal {G}]^{\bot}$. Thus (1) follows.

For $(2)$, note that
$
    B(I, Y)=B([I, \cdots, I], Y)=B([Y, I, \cdots, I], I)=0.
$
So $B|_{I\times I}$ is nondegenerate. This proves (2). The result (3) is trivial and (4) follows from the identity $[\mathcal {G}, \cdots, \mathcal {G}]^{\bot}=\mathcal {C}(\mathcal {G})$. ~$\Box$

\vspace{2mm}\noindent{\bf Lemma 2.4 } Let $W_{1}, W_{2}, \cdots,
W_{m}$ be subspaces of a metric $n$-Lie algebra $(\mathcal {G}, B)$.
Then $x \in [W_{1}, \cdots, W_{m}] ^{\bot}$ if and only if [$W_{1},
\cdots, W_{i-1}, x, W_{i+1}, \cdots, W_{m}] \subseteq W_{i} ^{\bot}$
for $1\le i \le m$.

\vspace{2mm}\noindent{\bf Proof }  Since $B(x, [W_{1}, W_{2},
\cdots,W_{m}]) =B(W_{i}, [W_{1}, \cdots, W_{i-1}, x, W_{i+1},
\cdots, W_{m}])$, we have that $x \in [W_{1}, \cdots, W_{m}]
^{\bot}$ if and only if $B(W_{i}, [W_{1}, \cdots, W_{i-1}, x,
W_{i+1}, \cdots, W_{m}])=0.$ The lemma follows.  ~$\Box$

\vspace{2mm}\noindent{\bf Theorem 2.5 } Let $(\mathcal {G}, B)$ be a
 $B$-irreducible metric $n$-Lie algebra. Then $\mathcal {G}$ is indecomposable.

\vspace{2mm}\noindent{\bf Proof } If $[\mathcal {G}, \cdots, \mathcal {G}]=0$, then $\dim \mathcal {G}=1$. The theorem holds. Suppose $[\mathcal {G}, \cdots, \mathcal {G}]\neq 0$. Let $\mathcal {G}=J\oplus H$ be a direct sum of ideals.  Then $[J, H, \mathcal {G}, \cdots, \mathcal {G}]=0$ and hence $[\mathcal {G},\cdots,\mathcal {G}]=[J, J, \mathcal {G}, \cdots, \mathcal
{G}]+[H, H, \mathcal {G}, \cdots, \mathcal {G}]$. So, without loss of generality, we may assume that $[J, J, \mathcal {G}, \cdots, \mathcal {G}]\neq 0$.

It follows from Lemma 2.3 that $B(J, [H, H, \mathcal {G}, \cdots, \mathcal
{G}])=0$ and $B(H, [J, J, \mathcal {G}, \cdots, \mathcal {G}])=0.$ Thus $J\subseteq [H, H, \mathcal {G}, \cdots, \mathcal {G}]^{\bot}$ and $H\subseteq [J, J, \mathcal {G}, \cdots, \mathcal {G}]^{\bot}$. Since $J^{\bot}\cap H^{\bot}=0,$ we obtain that $J^{\bot}
\cap [J, J, \mathcal {G}, \cdots, \mathcal {G}]=0$. Let $V$ be a
subspace of $J$ such that  $J=V\oplus (J^{\bot} \cap J)\oplus [J, J,
\mathcal {G}, \cdots, \mathcal {G}]$. Then $I=V\oplus [J, J, \mathcal {G}, \cdots, \mathcal {G}]$ is an ideal of $\mathcal {G}$ because $[J, H, \mathcal {G}, \cdots, \mathcal {G}]=0$.

We show that $I$ is nondegenerate. Let $x\in I$ and $B(x,I)=0$. Then $B(x, J)=0$, since $B(x, J \cap J^{\bot})=0$. Then  $x\in J^{\bot}\cap J\cap I$. It follows that $x=0.$ Therefore, $\mathcal {G}=I\subseteq J\subseteq
\mathcal {G}.$ Similarly, if $[H, H, \mathcal {G}, \cdots, \mathcal
{G}]\neq 0$, then $\mathcal {G}=H$.~~$\Box$

\vspace{2mm}\noindent{\bf Corollary 2.1 } Let $\mathcal {G}$ be an
$n$-Lie algebra with $B_1$ and $B_2$ being invariant scalar products
on $\mathcal {G}$. Then metric $n$-Lie algebra $(\mathcal {G}, B_1)$
is $B_1$-irreducible if and only if $(\mathcal {G}, B_2)$ is
$B_2$-irreducible.

\vspace{2mm}\noindent{\bf Proof } The result follows from Theorem 2.5.
~~$\Box$

\vspace{2mm}\noindent{\bf 3. The dimension of $B$-irreducible
metric $n$-Lie algebras }

To describe metric $n$-Lie algebras $(\mathcal {G}, B)$, we need
some notation and preliminary results. The centroid $\mathcal {G}$
is defined by
$$\Gamma (\mathcal {G})= \{\phi\in End(\mathcal {G})\mid
\phi ([x_{1}, x_{2}, \cdots, x_{n}]) = [\phi(x_{1}), x_{2}, \cdots,
x_{n}], \forall x_{i} \in \mathcal {G}\}. \eqno(3.1)
$$
It is easily seen that $\Gamma (\mathcal {G})$ is an associative algebra, and so it is a subalgebra of $gl(\mathcal {G})$. For more information on centroid see \cite{BAL}. Let
$$
\Gamma_{B}(\mathcal {G})=\{\phi \in \Gamma(\mathcal
{G}) \mid ~B(\phi(x),y)=B(x, \phi(y)), ~\mbox{for all}~  x, y\in
\mathcal {G}\}.
$$
Denote by $\Gamma^{0}_{B}(\mathcal {G})$ the subspace of
$\Gamma_{B}(\mathcal {G})$ spanned by $\{\phi \in
\Gamma_{B}(\mathcal {G}) \mid \phi ~\mbox{ is invertible}\}.$ We use
$\mathcal {F}(\mathcal {G})$ to represent the vector space spanned
by all invariant symmetric bilinear forms on $\mathcal {G}$ and
$B(\mathcal {G})$ the subspace of $\mathcal {F}(\mathcal {G})$
spanned by all nondegenerate invariant symmetric bilinear forms on
$\mathcal {G}$. The dimension of $B(\mathcal {G})$ is called the
{\it  metric dimension} of $\mathcal {G}$.

\vspace{2mm}\noindent{\bf Theorem 3.1 } Let $(\mathcal {G}, B)$ be a
metric $n$-Lie algebra. Then
\[
\dim \Gamma^{0}_{B}(\mathcal {G})=\dim\Gamma_{B}(\mathcal {G}) \quad\mbox{and}\quad
\dim B(\mathcal {G}) = \dim\mathcal {F}(\mathcal {G}).
\]

\vspace{2mm}\noindent{\bf Proof } Every  $\varphi \in \Gamma_{B}(\mathcal {G})$
corresponds to a matrix $M(\varphi)$  with respect to a basis. Note that $\varphi$ is nondegenerate if and
only if $\det {M}(\varphi)\neq 0$. Then the set of all nondegenerate elements
in $\Gamma^{0}_{B}(\mathcal {G})$ is a nonempty open subset in the topological
space $ \Gamma_{B}(\mathcal {G})$ with the usual topology. Therefore, $\dim \Gamma^{0}_{B}(\mathcal {G})=\dim\Gamma_{B}(\mathcal {G})$. By a similar discussion we get $\dim B(\mathcal {G}) = \dim\mathcal {F}(\mathcal
{G})$. ~$\Box$

\vspace{2mm}\noindent{\bf Lemma 3.2 } Let $(\mathcal {G}, B)$ be a
metric $n$-Lie algebra and $K$ be a bilinear form on $\mathcal {G}$.
Then there exists a unique linear transformation $D$ on $\mathcal G$
such that
$$
    K(x, y)=B(D x, y), ~\mbox{for all}~ x, y\in \mathcal {G}. \eqno(3.2)
$$
Furthermore, $K$ is symmetric and invariant if and only if $D$ satisfies
$$
    B( {D}(x), y)=B(x, {D}(y)), \eqno(3.3)
$$
and
$$
    D([x_{1}, x_{2}, \cdots, x_{n}])= [x_{1}, x_{2}, \cdots, D(x_{n})]\eqno(3.4)
$$
for all $x_{i}\in \mathcal {G}, 1\leq i\leq n$.

\vspace{2mm}\noindent{\bf Proof } Linear equations
$$
K(x_i, x_j)=B(D (x_i), x_j), ~\mbox{where}~ x_1, \dots, x_n ~\mbox{is a basis of}~ \mathcal {G}
$$
give rise to a linear transformation $D$ satisfying (3.2). Such a linear transformation is unique. It is routine to check that $K$ is symmetric if and only if $D$ satisfies (3.3). Since
$$
K(x_{n}, [x_{1}, \cdots, x_{n-1}, y]) = B(Dx_{n}, [x_{1},
\cdots, x_{n-1}, y]) = -B([x_{1}, x_{2}, \cdots, x_{n-1}, D(x_{n})],
y)
$$
and
$$
K([x_{1}, \cdots, x_{n-1}, x_n], y)=B(D[x_{1}, \cdots,
x_{n-1}, x_n], y]),
$$
we obtain that  $K$ is invariant if and only if (3.4) holds. ~$\Box$

\vspace{2mm}\noindent{\bf Theorem 3.3 } Let $(\mathcal {G}, B)$ be a
metric $n$-Lie algebra. Then we have $$\dim B(\mathcal {G})=\dim
\Gamma_{B}(\mathcal {G}).$$

\vspace{2mm}\noindent{\bf Proof } Let $B_{1}$ be an invariant scalar product on $\mathcal {G}$ and $B_{1}(x, y)=B(Dx, y)$ for some linear transformation $D$ on $\mathcal G$. Since $B_1$ is nondegenerate, by Lemma 3.2 we have $D\in \Gamma_{B}(\mathcal {G})$, it follows $\dim B(\mathcal {G}) \leq \dim \Gamma_{B}(\mathcal {G})$.

Conversely, let $\phi_1, \cdots, \phi_l$ be a basis of $\Gamma_{B}(\mathcal {G}).$ We define the invariant scalar products $B_{i}$ on $\mathcal G$ by
$$
B_{i}(x, y) = B (\phi_{i}(x), y), ~\forall x, y\in \mathcal {G}, ~\mbox{for}~  i=1, \cdots, l.
$$
Thus $\sum\limits_{i=1}^l k_{i} B_i(x,y)=B(\sum\limits_{i=1}^l k_{i}
\phi_{i}(x), y)$ for $k_1, \cdots, k_l\in \mathbb C$. It follows
that $\sum\limits_{i=1}^l k_{i} \phi_{i}=0$  if and only if
$\sum\limits_{i=1}^l k_{i} B_{i}=0$. Therefore, $\dim B(\mathcal
{G}) \geq \dim \Gamma_{B}(\mathcal {G})$. ~$\Box$

\vspace{2mm}\noindent{\bf Lemma 3.4 } Let $W$ be a subspace of
metric $n$-Lie algebra $(\mathcal {G}, B)$. Then

 \vspace{2mm}  (1)  $W$ is
an ideal of $\mathcal {G}$ if and only if $W^{\bot}$ is contained in
$\mathcal {C}_{\mathcal G}(W)$.

 \vspace{2mm} (2) All one-dimensional nondegenerate ideals are contained in
$\mathcal {C}(\mathcal G)$.

 \vspace{2mm}\noindent{\bf Proof }  If $W$ is an ideal of $\mathcal
{G}$, then $W^{\bot}$ is also an ideal. For every  $x\in W, y\in
W^{\bot}$,  $B([x, \mathcal {G}, \cdots, \mathcal {G}], y) =
B(\mathcal {G}, [x, y, \mathcal {G}, \cdots, \mathcal {G}])=0$.
Since $B$ is nondegenerate, $y\in \mathcal {C}_{\mathcal {G}}(W)$.
Conversely, if $W^{\bot}\subseteq \mathcal {C}_{\mathcal G}(W)$,
then $B(W^{\bot}, [W, \mathcal {G}, \mathcal {G}, \cdots, \mathcal
{G}])= B([W^{\bot}, W, \mathcal {G}, \cdots, \mathcal {G}], \mathcal
{G})=0$. It follows that $[W, \mathcal {G}, \mathcal {G}, \cdots,
\mathcal {G}]\subseteq (W^{\bot})^{\bot}=W$. This proves (1).

We prove $(2)$. If $I$ is a one dimensional nondegenerate ideal of $\mathcal {G}$,
then $[I, I, \mathcal {G}, \cdots, \mathcal {G}]=0$. By Lemma 2.3,
$$[I, \mathcal{G}, \cdots, \mathcal{G}]=[I, I\oplus I^{\bot},
\mathcal{G}, \cdots, \mathcal{G}]=[I, I, \mathcal{G}, \cdots,
\mathcal{G}] + [I, I^{\bot}, \mathcal{G}, \cdots, \mathcal{G}]=0.$$
This proves that  $I \subseteq \mathcal {C}(\mathcal G)$. ~$\Box$

\vspace{2mm}\noindent{\bf Theorem  3.5 } Let  $(\mathcal {G}, B)$ be
a metric $n$-Lie algebra and $\mathcal {G}={\mathcal S}\oplus
{\mathcal R}$ the Levi decomposition. Then $\mathcal {G}$ has no
strong semisimple ideals if and only if $\mathcal R$ is coisotropic
(i.e., $\mathcal R^{\bot} \subseteq \mathcal R$). Moreover, we have
$$
    \mathcal {C}_{\mathcal G}(\mathcal R)=\mathcal {C}(\mathcal{G}) \oplus {\mathcal R}^{\bot} \quad\mbox{and}\quad[\mathcal S, \cdots, \mathcal S, \mathcal R^{\bot}]=\mathcal R^{\bot}.
$$

\vspace{2mm}\noindent{\bf Proof } Suppose that $\mathcal G$ has no strong semisimple ideals. If $\mathcal S=0$, the result follows easily. If $\mathcal S\ne 0$, then the radical ${\mathcal R}_1$ of $\mathcal R^{\bot}$ is a characteristic ideal of $\mathcal R^{\bot}$, and hence $[{\mathcal R}_1, \mathcal S, \cdots, \mathcal S]\subseteq {\mathcal R}_1$.
Thus $\mathcal R^{\bot}$ has an $\mathcal \mathcal S$-module decomposition $\mathcal R^{\bot}={\mathcal R}_1\oplus S_1$, which is also a Levi decomposition of $\mathcal R^{\bot}$. If $S_1\ne 0$, it follows from Lemma 2.3 that $[S_1, \mathcal R, \mathcal {G}, \cdots, \mathcal {G}]=0$. In other words, $S_1$ is a strong semisimple ideal of $\mathcal {G}$, which is a contradiction. Therefore, $S_1=0$ and $\mathcal R^{\bot}={\mathcal R}_1\subseteq \mathcal R.$

Conversely, if $\mathcal R^{\bot}\subseteq \mathcal R,$ let $S_{1}$ be a strong
semisimple ideal of $\mathcal {G}$. By (2) in Lemma 2.3, $S_{1}$ is nondegenerate. Then $\mathcal {G}=S_{1}\oplus S_1^{\bot}$. Let $\mathcal {G}={\mathcal S}\oplus {\mathcal R}$ be the Levi decomposition of $\mathcal {G}$. Then $S_{1}\subseteq \mathcal S$ and $[S_{1}, \mathcal R, \mathcal {G}, \cdots, \mathcal {G}]\subseteq S_{1}\cap \mathcal R=0$. We have $B(\mathcal R,
S_{1})=B(\mathcal R, [S_{1}, \cdots, S_{1}])=B([S_{1}, \cdots, S_{1}, \mathcal R],
S_{1})=0$ since $[S_{1}, \cdots, S_{1}, \mathcal R] \subseteq S_{1}\cap \mathcal R=0$.
It follows that $S_{1}\subseteq \mathcal R^{\bot} \subseteq \mathcal R$. Therefore, $S_1=0.$ This shows that $\mathcal {G}$ has no strong semisimple ideals.

From Lemma 3.4, we have $\mathcal {C}_{\mathcal G}(\mathcal R)\supseteq \mathcal {C}(\mathcal{G})+ \mathcal R^{\bot}.$ Next we prove the inverse inclusion. Note that $\mathcal R^{\bot}$ is an $\mathcal S$-submodule. If
there is $x\in \mathcal R^{\bot}$ such that $[\mathcal S, \cdots, \mathcal S, x]=0$, then $B(x,
\mathcal {G})=B(x, \mathcal S)+B(x, \mathcal R)=B(x, \mathcal S)=B(x, [\mathcal S, \cdots, \mathcal S])=B(\mathcal S, [x,
\mathcal S, \cdots, \mathcal S])=0$. Therefore,  $x=0.$ It follows that $\mathcal R^{\bot}\cap
\mathcal {C}(\mathcal{G})=0$ and $[\mathcal S, \cdots, \mathcal S,
\mathcal R^{\bot}]=\mathcal R^{\bot}.$ Since
$$B(\mathcal R^{\bot}+
\mathcal {C}(\mathcal G), [\mathcal R, \mathcal {G}, \cdots,
\mathcal {G}])=B(\mathcal {C}(\mathcal G), [\mathcal R, \mathcal
{G}, \cdots, \mathcal {G}])=B(\mathcal R, [\mathcal {C} (\mathcal
G), \mathcal {G}, \cdots, \mathcal {G}])=0$$ and
 $\mathcal {C}_{\mathcal {G}}(\mathcal R)=[\mathcal R,
\mathcal {G}, \cdots, \mathcal {G}]^{\bot}$,  we get $\mathcal
R^{\bot} \oplus \mathcal {C}(\mathcal G) \subseteq \mathcal
{C}_{\mathcal {G}}(\mathcal R).$ Therefore, $\mathcal {C}_{\mathcal
{G}}(\mathcal R) = \mathcal R^{\bot} \oplus \mathcal {C}(\mathcal
G)$. ~$\Box$

\vspace{2mm}\noindent{\bf Theorem 3.6 } Let $(\mathcal {G}, B)$ be a
metric $n$-Lie algebra and $\mathcal {G}={\mathcal S}\oplus
{\mathcal R}$ be the Levi decomposition of $\mathcal {G}$. Then

(1) As $\mathcal S$-modules, we have
 $$\mathcal R^{\bot} \cong (\mathcal {G} / \mathcal R)^{*}
\cong \mathcal S^{*}\cong \mathcal S. \eqno(3.1)$$

(2) If $\mathcal S=S_1\oplus \cdots, \oplus  S_t$, where $S_i$ are simple ideals of $\mathcal S$ ($t\geq 1$), then there are $\mathcal S$-module homomorphisms $\tau_i: \mathcal S\rightarrow \mathcal R^{\bot}$ satisfying $\tau_i(S_j)=0$ ($i\neq j$), $\tau_i\mid_{S_i}$ is injective, and for all $x, y, x_1,
\cdots, x_n\in \mathcal S$
$$
B(\tau_i(x), y)=B(x, \tau_i(y)), \eqno(3.2)
$$
$$
\tau_i([x_1, \cdots, x_n])=[x_1, \cdots, x_{n-1}, \tau_i(x_n)]. \eqno(3.3)
$$

\vspace{2mm}\noindent{\bf Proof } If $\mathcal S=0,$ the result is trivial. If $\mathcal S\neq 0,$ it is clear that $(\mathcal {G}/\mathcal R)^{*} \cong \mathcal S^{*}\cong \mathcal S$ as $\mathcal \mathcal S$-modules since $\mathcal \mathcal S$ is strong semisimple. It suffices to show that $\mathcal R^\bot$ is isomorphic to $\mathcal S^*$ as $\mathcal S$-modules. Define a mapping $\phi$ of $\mathcal R^{\bot}$ to $\mathcal S^{*}$ by
$$
\phi(x)=B(x, -) \mid _{\mathcal S}, ~\mbox{ for all} ~x\in \mathcal R^{\bot}.\eqno(3.4)
$$
Since $B$ is nondegenerate, $\phi$ is surjective. If $\phi(x)=0$ for some
$x\in \mathcal R^{\bot}$, then $B(x, \mathcal S)=B(x, \mathcal S+\mathcal R)=B(x, \mathcal {G} )=0.$ It follows that $x=0,$ which shows that $\phi$ is injective. It is easily seen that $\phi$ is a linear isomorphism from vector space $\mathcal R^{\bot}$ to $\mathcal S^{\ast}$.
Moreover, for every $y_1, \cdots, y_{n-1}, z\in \mathcal S$ and $x\in
\mathcal R^{\bot}$ since
$$
(\mbox{ad}(y_1, \cdots, y_{n-1}).\phi(x))(z)  = -\phi(x)([y_1, \cdots, y_{n-1}, z])
$$
$$=-B(x, [y_1, \cdots, y_{n-1}, z])=B([y_1, \cdots, y_{n-1}, x],
z)=\phi([y_1, \cdots, y_{n-1}, x])(z),
$$ $\phi$ is an $\mathcal S$-module isomorphism. The result (1) follows.

Let $\mathcal S=S_1\oplus \cdots \oplus S_t$ and $\mathcal S^{\ast}=S_1^{\ast}\oplus \cdots \oplus S_t^{\ast},$  where $S_i$ are simple ideals of $\mathcal \mathcal S$ and $S_i^{\ast}$ are the dual space of $S_i$ (can also be seen as dual $\mathcal S$-module of $S_i$). From (1) we obtain a decomposition of irreducible $\mathcal S$-modules
$$
    \mathcal R^{\bot}={\mathcal R}_1\oplus \cdots \oplus {\mathcal R}_t,
$$
and the $\mathcal S$-module isomorphism $\pi_i: S_i\rightarrow \mathcal R_i$ satisfying
$$
    \pi_i([x_1, \cdots, x_n])=[x_1, \cdots, x_{n-1}, \pi_i(x_n)], ~ \mbox{ for all} ~x_1, \cdots, x_n\in S_i.\eqno(3.5)
$$
Thanks to (3.4) and Theorem 3.5, $ S_i\oplus {\mathcal R}_i$ is nondegenerate and
$$
    B(\pi_i(x), y)=B(x, \pi_i(y)), ~ \mbox{ for all} ~ x, y\in S_i.\eqno(3.6)
$$

Finally, extending $\pi_i$ to $\tau_i: \mathcal S\rightarrow \mathcal R^{\bot}$ by $\tau_i(x)=\pi_i(x)$ for $x\in S_i$ and $\tau_i(\sum\limits_{j\neq i}S_j)=0$ where $i=1, \cdots, t,$ we see $\tau_i$ as $\mathcal S$-module homomorphisms satisfy (3.2) and (3.3).
 ~$\Box$

\vspace{2mm}\noindent{\bf Lemma 3.7}$^{\cite{FO1}}$  An ideal $I$ is
maximal if and only if $\mathcal {G}/I$ is simple or
one-dimensional.

\vspace{2mm}\noindent{\bf Lemma 3.8 } If $\mathcal B$ is strong semisimple
or one dimensional, then every $n$-Lie algebra extension
$$
0\longrightarrow \mathcal A {\stackrel{p}{\longrightarrow}} \mathcal
{G}{\stackrel{\pi}{\longrightarrow}}  \mathcal B\longrightarrow 0
$$ is split, i.e., there is an $n$-Lie homomorphism $f:  \mathcal B\rightarrow
\mathcal {G}$ such that $\pi f=Id_{\mathcal B}.$

\vspace{2mm}\noindent{\bf Proof }  Let $\mathcal {G}=\mathcal R\oplus \mathcal S$
 be a Levi decomposition of $\mathcal {G}$ and $\mathcal S=S_{1} \oplus \cdots \oplus S_{m}$, where $S_i$ are simple subalgebras of $\mathcal {G}$. If $ \mathcal B$ is strong
semisimple, then the radical $\mathcal R$ of $\mathcal G$ is contained in the
kernel of $\pi$, and then  there is $k\leq m$ such that $\pi$ is an
isomorphism from $S_{i_1} \oplus \cdots \oplus S_{i_k}$ to $
\mathcal B.$ Let $\{y_1, \cdots, y_t\}$ be a basis of $S_{i_1} \oplus
\cdots \oplus S_{i_k}$. Then  $\{x_1 = \pi(y_1), \cdots, x_t = \pi(y_t)\}$
is a basis of $ \mathcal B$. Define a mapping $f:  \mathcal B\rightarrow \mathcal
{G}$ by $f(x_i)=y_i$ for $1\leq i\leq t$. Then $\pi f=Id_{
\mathcal B},$ which proves the Lemma.  If $\mathcal B$ is one-dimensional, the result is obvious.    ~ $\Box$

\vspace{2mm}\noindent{\bf Theorem 3.9 } Let $(\mathcal {G}, B)$ be a
metric $n$-Lie algebra  without strong semisimple ideals and
$\mathcal {G}={\mathcal S}\oplus {\mathcal R}$ the Levi
decomposition. Then a subalgebra $H$ of $\mathcal {G}$ is a minimal
ideal of $\mathcal {G}$ if and only if $H$ is an irreducible
$\mathcal S$-submodule of the centralizer $\mathcal
{C}_{\mathcal{G}}(\mathcal R)$.

\vspace{2mm}\noindent{\bf Proof } Since $[\mathcal {C}_{\mathcal
{G}}(\mathcal R), \mathcal R, \mathcal {G}, \cdots, \mathcal {G}]=0$, a subalgebra $H$
of $\mathcal {C}_{\mathcal {G}}(\mathcal R)$ is a minimal ideal if and only
if $H$ is an irreducible $\mathcal S$-module.

Now suppose a subalgebra $H$ of $\mathcal {G}$ is a minimal ideal. Then $H^{\bot}$ is a maximal ideal. Thanks to Lemma 3.7, $\mathcal
{G}/H^{\bot}$ is one-dimensional or simple. By
 Lemma 3.8 the short exact sequence $$0 \longrightarrow H^{\bot}
\longrightarrow \mathcal {G} \longrightarrow \mathcal {G} / H^{\bot}
\longrightarrow 0$$ is  split and there exists a subalgebra $S_{1}$
of $\mathcal {G}$ such that
$$\mathcal {G}=S_{1}\oplus H^{\bot}, \quad\mbox{as a direct sum of subalgebras}.\eqno(3.7)$$

If $\dim S_{1} >$1, then $S_{1}$ is simple and $\mathcal R\subseteq
H^{\bot}$. From Theorem 3.5 and Lemma 3.4,  we have $H\subseteq
\mathcal R^{\bot} \subseteq \mathcal {C}_{\mathcal {G}}(\mathcal R)$. Therefore, $H$ is an irreducible $\mathcal S$-module and contained in $\mathcal {C}_{\mathcal
{G}}(\mathcal R)$.

If $\dim S_{1}=1,$ then $[S_1, S_1, \mathcal {G}, \dots, \mathcal {G}]=0. $ By (3.7) we have
$$[H, \mathcal R, \mathcal {G}, \cdots, \mathcal {G}]=[H, \mathcal R, S_1, \mathcal {G}, \cdots, \mathcal {G}].
$$
We claim that the right hand is zero. Suppose not, as $H$ is a minimal ideal we see $[H, \mathcal R, S_1, \mathcal {G}, \cdots,
\mathcal {G}]=H$. Then $B(H, \mathcal {G})=B(H, S_1)=B([H, \mathcal R, S_1,
\mathcal {G}, \cdots, \mathcal {G}], S_1)=B(H, [S_1, \mathcal R, S_1,
\mathcal {G}, \cdots, \mathcal {G}])=0,$ which contradicts that $B$ is nondegenerate on $\mathcal G$. Therefore, $H\subseteq \mathcal {C}_{\mathcal {G}}(\mathcal R)$ and
 $H$ is an irreducible
$\mathcal S$-module. ~$\Box$

\vspace{2mm}\noindent{\bf Lemma 3.10 } Let $(\mathcal {G}, B)$ be a
metric $n$-Lie algebra. Then $\mathcal {C}(\mathcal G)\subseteq
\mathcal {C}(\mathcal G)^{\bot}$ if and only if $\mathcal {G}$ has
no one-dimensional nondegenerate ideals.

\vspace{2mm}\noindent{\bf Proof } If   $\mathcal {C} (\mathcal G)$ is not contained in $\mathcal {C}(\mathcal G)^{\bot}$, then there exists $x\in\mathcal {C}(\mathcal G)$ such that $B(x, x)\neq 0.$ We obtain a one dimensional nondegenerate ideal $I=\mathbb C x$.

Conversely, suppose that $\mathcal {C}(\mathcal G)\subseteq \mathcal {C}(\mathcal G)^{\bot}$ and that $J=\mathbb C x$ is a one dimensional ideal of $\mathcal G$. Then we have $[\mathcal G, \cdots, \mathcal G, x]=0$ or $[\mathcal G, \cdots, \mathcal G, x]=\mathbb Cx.$ If $[\mathcal G, \cdots, \mathcal G, x]=0$, then
$x\in C(\mathcal G)\subseteq \mathcal {C}(\mathcal G)^{\bot}$, it follows $B(J,
J)=0.$ If $[\mathcal G, \cdots, \mathcal G, x]=\mathbb Cx,$ Then $B(x, x)\subseteq B( x, [\mathcal G, \cdots, \mathcal G, x])=B(\mathcal G, [x, x, \mathcal G, \cdots, \mathcal G])=0. $  Therefore, $J$ is
degenerate. ~$\Box$

\vspace{2mm}\noindent{\bf Theorem 3.11 } Let $(\mathcal {G}, B)$ be
a $B$-irreducible metric $n$-Lie algebra and $m(\mathcal {G})$ the
number of minimal ideals in a decomposition of ${\rm Soc}(\mathcal
G)$. Then $\dim B(\mathcal {G})\geq m(\mathcal {G})+1.$

\vspace{2mm}\noindent{\bf Proof }  Let $\mathcal {G}={\mathcal S}\oplus {\mathcal R}$ be
the Levi decomposition of $\mathcal {G}$. If $\mathcal S=0$ then $m(\mathcal {G})=\dim \mathcal {C}(\mathcal {G})$ by Theorem 3.9. The result follows from Theorem 3.5.

If $\mathcal S\neq 0$, then from Theorem 3.5, Lemma 3.10 and invariant property of $B$, we have $\mathcal R^{\bot}\subseteq \mathcal R=(\mathcal R^{\bot})^{\bot}$ and $\mathcal {C} (\mathcal
{G})\subseteq (\mathcal {C} (\mathcal {G}) )^{\bot}$. It follows that
$\mathcal {C}_{\mathcal {G}}(\mathcal R)\subseteq \mathcal {C}_{\mathcal
{G}}(\mathcal R)^{\bot}$,
$$
    B(\mathcal R^{\bot},\mathcal R^{\bot})=B(\mathcal {C}(\mathcal {G}), \mathcal {C}(\mathcal
    {G}))=B(\mathcal {C}(\mathcal {G}), \mathcal R^{\bot})=B(S, \mathcal C(\mathcal {G}))=0.
$$
If $\mathcal {C} (\mathcal {G})\neq 0$, then there is a subspace $V\subseteq \mathcal R$ such that $\mathcal {C} (\mathcal {G})\oplus V$ is nondegenerate (as a direct sum of subspaces). Thus $\mathcal {G}=V\oplus [\mathcal{G}, \cdots, \mathcal {G}]$ again by Theorem 3.5 and Lemma 3.10.

Let $\mathcal S=S_1\oplus \cdots \oplus S_t$, where $S_i$ are simple
ideals of $\mathcal S$.  By Theorem 3.6, $\mathcal R^{\bot}$ has an
irreducible $\mathcal S$-module decomposition $\mathcal
R^{\bot}={\mathcal R}_1\oplus \cdots \oplus {\mathcal R}_t$. From
Theorem 3.9, $$m(\mathcal {G})=t+\dim \mathcal {C} (\mathcal
{G}).\eqno(3.8)$$

Extend $\tau_i$ in Theorem 3.6 to a mapping $\phi_i$ of $\mathcal {G}$ to itself by $\phi_i(x)=\tau_i(x)$  for $x\in S$ and $\phi_i(x)=0$ for $x\in \mathcal R$. We obtain, for all  $x, y, x_1, \cdots, x_n\in \mathcal {G}$,
$$
B(\phi_i(x), y)=B(x, \phi_i(y)), \quad \phi_i([x_1, \cdots, x_n])=[x_1,
\cdots, x_{n-1}, \phi_i(x_n)].
$$
In other words, $\phi_i\in \Gamma_B(\mathcal {G})$ for $1\leq i\leq t.$ Let $\{y_1, \cdots, y_s\}$ be a basis of $\mathcal {C}(\mathcal {G})$. Then there is a basis $\{x_1, \cdots, x_s\}$ of $V$ such that $B(x_i, y_j)=\delta_{ij}.$  Define $\psi_j: \mathcal {G}\rightarrow \mathcal {G}$ by
$$
    \psi_j(x_i)=\delta_{ij} y_j\quad\mbox{and}\quad \psi_j([\mathcal {G}, \cdots, \mathcal {G}])=0, \quad\mbox{for } i=1, \cdots, s.
$$
So $\psi_j(\mathcal {G})\subseteq \mathcal {C}(\mathcal {G})$
and  $B(\psi_j(x), y)=B(x, \psi_j(y))$ for $x, y\in \mathcal
{G}.$ It follows that $\psi_j\in \Gamma_B(\mathcal {G})$.

In sum, $\{\phi_i, \psi_j \in \Gamma_{B}(\mathcal {G}) \mid 1\leq i\leq t, 1\leq j\leq s\}$ is linear independent. Since the identity map $Id_{\mathcal {G}}$ is contained in $\Gamma_{B}(\mathcal {G})$, we have $\dim B(\mathcal {G})\geq m(\mathcal {G})+1$ by (3.8).   ~~$\Box$

 By the above discussions, we  obtain following results directly.

\vspace{2mm}\noindent{\bf Corollary 3.1 } Let $(\mathcal {G}, B)$ be
a metric $n$-Lie algebra with dim $\mathcal {G}>$1. Then

(1) $\dim B(\mathcal {G})\geq \dim\mathcal {C}(\mathcal {G})+1$.

(2) If $\mathcal {G}= \oplus_{i=1}^p ~\mathcal {G}_{i}$ is an
orthogonal direct sum  of $B$-irreducible ideals, then $\dim
B(\mathcal {G})\geq \sum_{i=1}^p \dim B(\mathcal {G}_{i})$.
Moreover, $\dim B(\mathcal {G})\geq p + \sum_{i=1}^p m(\mathcal
{G}_{i})$.

(3) If $\dim B(\mathcal {G})=2,$ then $\mathcal {G}$ is a local $n$-Lie algebra if and only if $\mathcal {G}$ is either a solvable $n$-Lie algebra whose center is one-dimensional or a perfect $n$-Lie algebra with a simple Levi factor.

\vspace{2mm}\noindent{\bf Proof } The results (1) and (2) follow from the
Theorem 3.11. We prove (3) below. Let $\mathcal {G}$ be a local $n$-Lie
algebra. Then $m(\mathcal {G})=1.$ If $\mathcal {G}$ is solvable, by Theorem 3.9 the only minimal ideal of $\mathcal {G}$ is contained in $\mathcal C(\mathcal G)$. Therefore, $\dim \mathcal {C}(\mathcal {G})=1.$  If $\mathcal {G}$ is not solvable. Let $\mathcal G={\mathcal S}\oplus {\mathcal R}$ be the Levi decomposition and $\mathcal S=S_1\oplus \cdots\oplus S_t$, where $S_i$ are simple subalgebras ($1\le i \le t)$. Since $\dim B(\mathcal
{G})=2$, by (3.8) and Theorem 11, we have $m(\mathcal G)=t+\dim \mathcal {C}(\mathcal {G})=1$. Thus  $\dim \mathcal {C}(\mathcal {G})=0$ and $t=1$. It follows from Theorem 3.5 that $\mathcal R=\mathcal R^{\bot}$ and $[\mathcal {G}, \cdots, \mathcal G]=[S, \cdots, S]\oplus [S, \cdots, S,\mathcal R^{\bot}]={\mathcal S}\oplus {\mathcal R}^{\bot}=\mathcal G$, which shows that $\mathcal G$ is a perfect $n$-Lie algebra with a simple Levi factor. The converse follows from Theorem 3.11. ~$\Box$

\vspace{2mm}\noindent{\bf 4.  Minimal ideals of metric $n$-Lie
algebras}

In this section we study structures of metric $n$-Lie algebras by
means of minimal ideals.

\vspace{2mm}\noindent{\bf Lemma 4.1 } Let $\mathcal {G}$ be an $n$-Lie algebra and $J_{k}$ ideals of $\mathcal {G}$ for $k=1, \cdots, m$. Suppose that $\mathcal {G}=\oplus_{k=1}^{m} J_{k}$ is a direct sum of vector spaces. If $J$ is an ideal satisfying $[J, \mathcal {G}, \cdots, \mathcal {G}]=J$ or $\mathcal {C}(\mathcal {G} / J)= \{0\}$, then $J=\oplus_{k=1}^{m}(J_{k}\cap J)$.

\vspace{2mm}\noindent{\bf Proof }  If $[J, \mathcal {G}, \cdots,\mathcal {G}]$ $=J$, then  $J = \oplus_{k=1}^{m} [J, J_{k}, \cdots,J_{k}]\subseteq \oplus_{k=1}^{m} (J_{k}\cap J)$. Therefore, $J=\oplus_{k=1}^{m} (J_{k}\cap J)$. If $\mathcal {C}(\mathcal {G} / J)=\{0\}$. Let $\varphi :\mathcal {G} \rightarrow \mathcal {G}/ J$ be the canonical mapping and $x=\sum_{k=1}^{m} x_{k}\in J$, where $x_{k}\in J_{k}, 1\leq k\leq m$.  Since $[x_{k}, \mathcal {G}, \cdots, \mathcal {G}]=[x_{k}, J_{k}, \cdots, J_{k}]$ and $ [x_{k}, J_{j}, \mathcal {G}, \cdots, \mathcal {G}]=0$ if $ k\neq j$, we have $[x_{k}, \mathcal {G}, \cdots, \mathcal {G}]=[x, J_{k}, \cdots, J_{k}]\subseteq J.$ Then
$$
    [\varphi(x_k), \varphi(\mathcal {G}), \cdots, \varphi(\mathcal {G})]=\varphi([x_k, \mathcal {G}, \cdots, \mathcal {G}])=0,
$$
which shows that $\varphi(x_k)\in \mathcal {C}(\mathcal {G} / J).$ Thus $x_k\in J$. This completes the proof. ~~$\Box$

From \cite{BM}, an $n$-Lie algebra $\mathcal {G}$ is strong semisimple if
and only if $\mathcal {G}$ can be decomposed into a direct sum of
simple ideals. The following lemma provides a property of the
maximal strong semisimple ideals in metric $n$-Lie algebras.

\vspace{2mm}\noindent{\bf Lemma 4.2 } Let $\mathcal {G}$ be a metric
$n$-Lie algebra with simple ideals. Then $\mathcal {G}$ has a unique
maximal strong semisimple ideal $\mathcal {S}(\mathcal{G})$ and
$$
    \mathcal {G}=\mathcal {S}(\mathcal {G})\oplus \mathcal {S}(\mathcal
    {G})^{\bot}.
$$

\vspace{2mm}\noindent{\bf Proof } Suppose $\mathcal M$ is the set of
all simple ideals of $\mathcal G$ and $\mathcal {S}(\mathcal
{G})=\sum_{J\in \mathcal M}J$. For every $J_1, J_2 \in \mathcal M$,
since $[J_1, J_2, \mathcal G, \cdots, \mathcal G]\subseteq J_1\cap
J_2 ,$ we have  $J_1=J_2$ or  $[J_1, J_2, \mathcal G, \cdots,
\mathcal G]=0$. It follows that $\mathcal {S}(\mathcal {G})$ is the
direct sum of all simple ideals of $\mathcal G$. Therefore,
$\mathcal {S}(\mathcal {G})$ is the unique maximal strong semisimple
ideal.

Let $\mathcal {G}={\mathcal S}\oplus {\mathcal R}$ be the Levi decomposition of
$\mathcal {G}$. Then $\mathcal {S}(\mathcal {G})$ is contained in $\mathcal S$, and $\mathcal S$ has a  decomposition  $\mathcal S=\mathcal {S}(\mathcal
{G})\oplus S_1$, where $S_1$ is an ideal of $\mathcal S$. Let  $\mathcal
{G}_1=S_1\oplus {\mathcal R}$. Then $\mathcal R$ is the radical of $\mathcal
{G}_1$, and
$$
[\mathcal {S}(\mathcal {G}), \mathcal R, \mathcal {G}, \cdots, \mathcal
{G}]=0, \quad [\mathcal {S}(\mathcal {G}), S_1, \mathcal {G}, \cdots,
\mathcal {G}]=0.
$$
It follows that
$$
 [\mathcal {G}_1, \mathcal {G}, \cdots, \mathcal
{G}]=[\mathcal {G}_1, \mathcal {G}_1, \mathcal {G}, \cdots, \mathcal
{G}]=S_1+[ \mathcal R, \mathcal {G}, \cdots, \mathcal {G}]\subseteq \mathcal
{G}_1,
$$
which shows that  $\mathcal {G}_1$ is an ideal of $\mathcal {G}$ and has no strong semisimple ideals. By Lemma 2.3, we see that $\mathcal {S}(\mathcal {G})$ is nondegenerate and $\mathcal
{G}_1=\mathcal {S}(\mathcal {G})^{\bot}$. ~~$\Box$

\vspace{2mm}\noindent{\bf Lemma 4.3 } Let $(\mathcal {G}, B)$ be a
nonsimple  metric $n$-Lie algebra with $\dim \mathcal {G} > 1$ and
$J$ be a non-zero minimal ideal of $\mathcal {G}$. Then

(1) $J$ is either an abelian or a simple ideal of $\mathcal {G}$.

(2) $J\subseteq [\mathcal {G}, \cdots, \mathcal {G}]$ or $J\subseteq
\mathcal {C}(\mathcal {G})$.

(3) If $H\neq J$ is a minimal ideal, then  $[J,H,\mathcal {G}, \cdots, \mathcal {G}]=0$.

\vspace{2mm}\noindent{\bf Proof } (1) From the minimality of $J$ we
obtain that $J\cap J^{\bot}=J$ or $J\cap J^{\bot}=0$. If $J\cap
J^{\bot}=J$, then $J\subseteq J^{\bot}$. It follows from Lemma 2.3
that $[J, J, \mathcal {G}, \cdots, \mathcal {G}] \subseteq [J,
J^{\bot}, \mathcal {G}, \cdots, \mathcal {G}]=0$. Then $J$ is an
abelian ideal. If $J\cap J^{\bot}=0,$ then $\mathcal {G}=J\oplus
J^{\bot}.$ Note that a subspace $I\subseteq J$ is an ideal of $J$ if
and only if $I$ is an ideal of $\mathcal {G}$, since $[I, \mathcal
{G}, \cdots, \mathcal {G}]=[I, J, \cdots, J]$. Therefore, if $\dim
J=1$, then $J$ is abelian; if $\dim J >1$, then $J$ is simple.

The result (2)  can be seen from $[J, \mathcal {G}, \cdots, \mathcal {G}]=J$ or $[J, \mathcal {G}, \cdots, \mathcal {G}]=0$. The result (3) is trivial. ~~$\Box$

\vspace{2mm}\noindent{\bf Theorem 4.4 } Let $(\mathcal {G}, B)$ be a
nonsimple metric $n$-Lie algebra with $\dim \mathcal {G} > 1$ and
$J$ a minimal ideal of $\mathcal {G}$. Then $\mathcal {G} /J^{\bot
}$ is a one-dimensional metric $n$-Lie algebra if and only if
$J\subseteq \mathcal {C}(\mathcal {G})$ and $\dim J=1$.

\vspace{2mm}\noindent{\bf Proof } If $\mathcal {G} /J^{\bot}$ is a one-dimensional metric $n$-Lie algebra, then there is a nonzero vector $x\in \mathcal {G}$ such that $\mathcal {G}=J^{\bot}\oplus \mathbb C x$ as a direct sum of vector spaces. In the case that $J$ is nondegenerate, we have $\mathcal {G}=J\oplus J^{\bot}$. Thus $\dim J=1$. By Lemma 3.4, $[J, \mathcal {G}, \cdots, \mathcal {G}]$ $= [J, J+J^{\bot}, \mathcal {G}, \cdots, \mathcal {G}]=0$, that is, $J\subseteq \mathcal {C}(\mathcal {G})$. If $J$ is degenerate, then $J\cap J^{\bot} \neq 0$ and $J\subseteq J^{\bot}.$ So the subspace  $J\oplus \mathbb C x$ is nondegenerate. It follows that $\dim J=1$, and $J\subseteq \mathcal {C}(\mathcal {G})$.

Now, suppose that $J=\mathbb C  z \subseteq \mathcal {C}(\mathcal {G}).$ If $J$ is nondegenerate, we have $\mathcal {G}=J\oplus J^{\bot}$; and $\mathcal {G}/J^{\bot}$ is isomorphic to $J$ as one-dimensional metric $n$-Lie algebras. If $J$ is degenerate, then $J\subseteq J^{\bot}.$ Since $J$ is a minimal ideal, $J^{\bot}$ is maximal. Then  $\mathcal {G}/J^{\bot}$
is one-dimensional or simple $n$-Lie algebra. By Lemma 3.8, there is a subalgebra $S_1$ of $\mathcal {G}$ being one-dimensional or simple such that $\mathcal {G}=J^{\bot}\oplus S_1$ as a direct sum of subalgebras. We affirm that $S_1$ is not simple. Otherwise, $S_1=[S_1, \cdots, S_1]$ and hence $B(J, S_1)=0$ and $B(J, \,\mathcal {G})=0,$ which contradicts to that $B$ is nondegenerate. It follows that $S_1= \mathbb C y$ and $J\oplus \mathbb C y$ is nondegenerate. Thus, if $\mathbb C y$ is nondegenerate, the result follows. If $\mathbb C y$ is degenerate, then $\mathcal {G}/J^{\bot}$ is isomorphic to $\mathbb C (z+y)$ as a metric  one-dimensional $n$-Lie algebra.
 ~~$\Box$

\vspace{2mm}\noindent{\bf Theorem 4.5 } Let $(\mathcal {G}, B)$ be a
$B$-irreducible nonsimple metric $n$-Lie algebra with $\dim \mathcal
{G} >1$ and $J$ a nonzero minimal ideal of $\mathcal {G}$. Then
$J\subseteq \mathcal {C}_{\mathcal {G}}(\mathcal R)\cap [\mathcal
{G}, \cdots, \mathcal {G}]$.

\vspace{2mm}\noindent{\bf Proof } Since  $\mathcal {G}$ is a
$B$-irreducible nonsimple $n$-Lie algebra, by Theorem 3.9 we have
$J\subseteq \mathcal {C}_{\mathcal {G}}(\mathcal R)$. Note that
$[J, \mathcal {G}, \cdots, \mathcal {G}]$ is an ideal contained in
$J$. We have $[J, \mathcal {G}, \cdots, \mathcal {G}]=0$ or $[J,
\mathcal {G}, \cdots, \mathcal  {G}]=J$. If $[J, \mathcal {G},
\cdots, \mathcal {G}]=J$, then the  result follows.  If  $[J,
\mathcal {G}, \cdots, \mathcal {G}]=0$, then $J\subseteq \mathcal
{C}(\mathcal {G})$. Since $\mathcal {G}$ is $B$-irreducible,  by
Lemma 3.10, $\mathcal {C}(\mathcal {G})$ is  isotropic. Therefore,
$J\subseteq \mathcal {C}(\mathcal {G}) \subseteq \mathcal
{C}(\mathcal {G})^{\bot}=[\mathcal {G}, \cdots, \mathcal {G}]$,
which completes the proof.     ~~$\Box$

\vspace{1mm} A direct calculation yields the following lemma. We omit its details.

\vspace{1mm}\noindent{\bf Lemma 4.6 } Let $\mathcal {G}$ be a metric
$n$-Lie algebra and $\mathcal {S}(\mathcal {G})$ be the maximal
strong semisimple ideal. Suppose that $\mathcal {S}(\mathcal
{G})^\bot = S_1\oplus \mathcal {R}_1$ is the Levi decomposition of
$\mathcal {S}(\mathcal {G})^\bot$. Then we have

 (1) Every ideal of $\mathcal {S}(\mathcal {G})^\bot$ is an ideal of $\mathcal {G}$, and $\mathcal {R}_1 = \mathcal {R}$, the radical of $\mathcal G$.

 (2) $\mathcal {S}(\mathcal {G})^\bot$ has not strong semisimple ideals.

(3) Let $\mathcal {A}(\mathcal {G})$ be the sum of all minimal abelian ideals of $\mathcal {G}$. Then  $\mathcal {A}(\mathcal {G})=Soc(\mathcal {S}(\mathcal {G})^\bot)$.

(4) $\mathcal {C}(\mathcal {G}) = \mathcal {C}(\mathcal {S}(\mathcal {G})^\bot)$.

 (5) $\mathcal {C}_{\mathcal G}(\mathcal {R}) =\mathcal {S}(\mathcal {G})\oplus \mathcal {C}_{\mathcal {S}(\mathcal {G})^\bot}(\mathcal {R})$.

(6) $\mathcal {R}^{\bot}=\mathcal {S}(\mathcal {G}) \oplus W$, where
$W$ is the orthogonal complement of $\mathcal {R}$ in $\mathcal
{S}(\mathcal {G})^\bot$.

\vspace{2mm}\noindent{\bf Theorem 4.7 } Let $(\mathcal {G}, B)$ be a
metric $n$-Lie algebra  and $\mathcal {G}={\mathcal S}\oplus
{\mathcal R}$ be its Levi decomposition. Then
$$
    {\rm Soc}(\mathcal {G})= \mathcal {C}_{\mathcal {G}}(\mathcal R)=
    \mathcal R^{\bot}\oplus\mathcal {C}(\mathcal {G}).
$$

\vspace{2mm}\noindent{\bf Proof } Let $\mathcal {S}(\mathcal {G})$
be the maximal strong semisimple ideal and $\mathcal {A}$ be the set
of all nonsimple minimal ideals of $\mathcal {G}$ and $\mathcal
{A}(\mathcal {G})=\sum _{J\in \mathcal {A}} J$. By Lemma 4.2 and
Lemma 4.3, we have $\mathcal {A}(\mathcal {G})$ is an abelian ideal
and
$$
Soc(\mathcal {G})=\mathcal {S}(\mathcal {G})\oplus \mathcal
{A}(\mathcal {G}), ~~~\mathcal {A}(\mathcal {G})\subseteq \mathcal R, ~~~
\mathcal {S}(\mathcal {G}) \subseteq \mathcal {C}_{\mathcal
{G}}(\mathcal R).
$$
Then $\mathcal {A}(\mathcal {G})\subseteq \mathcal {C}_{\mathcal {S}(\mathcal {G})^\bot}(\mathcal R)\subseteq \mathcal {C}_{\mathcal {G}}(\mathcal R)$ by Theorem 3.9 and (3) in Lemma 4.6. Therefore,
$$
    {\rm Soc}(\mathcal {G})=\mathcal {S}(\mathcal {G})\oplus \mathcal {A}(\mathcal {G}) \subseteq \mathcal {C}_{\mathcal {G}}(\mathcal R).
$$
Let $\mathcal {S}(\mathcal {G})^\bot = S_1\oplus \mathcal {R}$ be the Levi decomposition of $\mathcal {S}(\mathcal {G})^\bot$ and $W$ the orthogonal complement of $\mathcal {R}$ in $\mathcal {S}(\mathcal {G})^\bot$.  It follows from Theorems 3.5 and 3.6 that $\mathcal {C}_{\mathcal {S}(\mathcal {G})^\bot}(\mathcal R) =\mathcal {C}(\mathcal {S}(\mathcal {G})^\bot) \oplus W = \mathcal {C}(\mathcal {G})\oplus W$ and that $W$ has an irreducible $S_1$-submodule decomposition $W = W_1 \oplus \cdots\oplus W_t$. Again by
Theorem 3.9 we know that $W_i$ are minimal ideals of $\mathcal {S}(\mathcal {G})^\bot$. Thus $\mathcal {C}_{\mathcal {S}(\mathcal {G})^\bot}(\mathcal R)\subseteq Soc(\mathcal {G})$. Therefore, $$\mathcal {C}_{\mathcal {G}}(\mathcal R) \subseteq \mathcal {S}(\mathcal {G})\oplus \mathcal {A}(\mathcal {G}).
$$
This completes the proof.   ~~$\Box$

\end{document}